\documentclass{article}
\usepackage{amssymb,amsmath,amsthm}
\usepackage{amsfonts}
\usepackage{amscd}
\newtheorem{thm}{Theorem}
\begin{document}

\title{About estimations of difference for the partial
integro-differential equation with small parameter at the leading
derivative}
\author{I.~Kopshaev
\thanks{Institute of mathematics of NAS of KAZAKHSTAN, 125 Pushkina str.,
050010 Almaty, KAZAKHSTAN email: {\em{kopshaev@math.kz}}}}

\maketitle
\date

\begin{abstract}
This paper is devoted to the study of the singularly perturbed
second order partial integro-differential equations. The
estimation of the solutions of Cauchy problem is obtained.
\end{abstract}

\section{Introduction}

Constructing of asymptotic decompositions of solutions of
singularly perturbed differential equations has great theoretical
and practical importance. In this line of investigations, the
fundamental results obtained by A.N. Tihonov \cite{Tihonov}, A.B.
Vasyliyeva \cite{Vasyl}, K.A. Kasymov \cite{Kasymov}, M.I.
Imanaliyev \cite{Imanaly}, L.A. Lyusternik \cite{Vishik} and
others.

However, the results, obtained in \cite{Tihonov}-\cite{Vishik},
still not generalized for the partial second order
integro-differential equations of Volterra type.

The aim of this paper is the study of the Cauchy problem with
initial jump for the singularly perturbed partial second order
integro-differential equations.

\section{Preliminaries}
Consider in $G=\{(t,x): 0 \leq t \leq 1, \lambda \leq x \leq
\lambda + 1\}$ following problem
$$
L_{\varepsilon}y=\varepsilon H^{2}[y]+A(t,x)H[y]+B(t,x)y=
$$
\begin{equation}
\label{kop_1} =F(t,x)+
\int_{0}^{t}\left(K_{1}(t,s,x)H[y(s,x)]+K_0(t,s,x)y(s,x)\right)ds,
\end{equation}
\begin{equation}
\label{kop_2} y(0,x,\varepsilon)=\pi_{0}(x), \qquad \varepsilon
\cdot  y_{t}^{'}(0,x,\varepsilon)=\pi_{1}(x).
\end{equation}
Here $\varepsilon>0$ - a small parameter, $t,x$ - independent
variables, $y=y(t,x,\varepsilon)$ - unknown function, $A(t,x)$,
$B(t,x)$, $F(t,x)$, $K_{i}(t,s,x)$ and $\pi_{i}, (i=0,1)$ -
functions given in G, operators
$$H[y]=<e(t,x) \cdot grad y>, \quad H^{2}[y]=H[H[y]],$$
where $<\cdot>$ denotes inner product of vectors $e(t,x) = (1,
Q(t,x))$ and $grad y=\left(\frac{\partial y}{\partial t},
\frac{\partial y}{\partial x}\right)$, $Q(t,x)$ is also given in
G, the function $\lambda(t)$ is a solution of characteristic
equation
\begin{equation}
\label{kop_3} \frac{dx}{dt}=Q(t,x).
\end{equation}
Consider also disturbed problem
$$
L_{0}y_0=A(t,x)H[y_0]+B(t,x)y_0=
$$
\begin{equation}
\label{kop_4} =F(t,x)+
\int_{0}^{t}\left(K_{1}(t,s,x)H[y_0(s,x)]+K_0(t,s,x)y_0(s,x)\right)ds,
\end{equation}
\begin{equation}
\label{kop_5} y_0(0,x)=\pi_{0}(x).
\end{equation}
obtained from (\ref{kop_1}), (\ref{kop_2}) when $\varepsilon=0$.

Suppose, that

1) $A(t,x)$, $B(t,x)$, $F(t,x)$, $K_{i}(t,s,x)$, $\lambda(t)$ and
$\pi_{i}, (i=0,1)$ - continuous functions in G.

2) conditions
$$
\inf\limits_{(t,x) \in G} A(t,x) \ge \gamma > 0, \quad
\inf\limits_{(t,x) \in G} Q(t,x) \ge \sigma > 0, \quad
\inf\limits_{(t,x) \in G} \pi_{i}(x) \ge \sigma > 0,
$$
$$
\lambda(0) = 0, \lambda(1) = 1,
$$
is satisfied, where $\gamma$ and $\sigma$ - some real numbers.

\section{Estimation of difference}
In this section, I prove estimations of the difference between
perturbed and unperturbed Cauchy problems.

Using theorem 2 from \cite{Tazhim}, is not difficult to prove,
that solution $y(t,x,\varepsilon)$ of problem (\ref{kop_1}),
(\ref{kop_2}) doesn't goes to solution of problem (\ref{kop_4}),
(\ref{kop_5}) when $\varepsilon \to 0$.

Consider the following problem
\begin{equation}
\label{kop_6} L_{0}y_0=F(t,x)+
\int_{0}^{t}\left(K_{1}(t,s,x)H[y_0(s,x)]+K_0(t,s,x)y_0(s,x)\right)ds
+ \Delta(t,x),
\end{equation}
\begin{equation}
\label{kop_7} y_0(0,x)=\pi_{0}(x)+ \Delta_0(x),
\end{equation}
where $\Delta(t,x)$, $\Delta_0(x)$ - not for a while yet unknown
functions.

Function $\Delta_0(x)$ is call to be named initial jump of
solution of problem (\ref{kop_1}), (\ref{kop_2}), function
$\Delta(t,x)$ - initial jump of integral term of equation
(\ref{kop_1}).

Suppose, that solution $y_0(t,x)$  of problem (\ref{kop_6}),
(\ref{kop_7}) when $t=t_0=\frac{\varepsilon}{\gamma} |\ln
\varepsilon|$ satisfies the condition
\begin{equation}
\label{kop_8} y_0(t_0,x)=y(t_0,x,\varepsilon), \quad x \in G.
\end{equation}

\begin{thm}
Let conditions 1),2) be satisfied. Then for the difference between
solution of problem (\ref{kop_1}), (\ref{kop_2}) and the solution
of problem (\ref{kop_6}), (\ref{kop_7}) in $G_1 \subset G$ have
place following estimations
$$
|y(t,x,\varepsilon)-y_0(t,x)| \leq K \cdot \varepsilon \cdot |\ln
\varepsilon|+K \cdot \varepsilon \cdot
|y_{t}^{'}(t_0,x,\varepsilon)|+K \cdot \max\limits_{(t,x \in G)}
|K_{1}(t,0,x) \cdot \Delta_{0}(x) - \Delta(t,x)|,
$$
$$
|y_{t}^{'}(t,x,\varepsilon)-y_{0t}^{'}(t,x)| \leq K \cdot
\varepsilon \cdot |\ln \varepsilon|+K \cdot \varepsilon \cdot
|y_{t}^{'}(t_0,x,\varepsilon)|+
$$
$$
+ K \cdot \max\limits_{(t,x \in G)} |K_{1}(t,0,x) \cdot
\Delta_{0}(x) - \Delta(t,x)| + K \cdot
\left(1+|y_{t}^{'}(t_0,x,\varepsilon)| \right) \cdot
e^{-\frac{\gamma}{\varepsilon}(t-t_0)},
$$
$$
|y_{x}^{'}(t,x,\varepsilon)-y_{0x}^{'}(t,x)| \leq K \cdot
\varepsilon \cdot |\ln \varepsilon|+K \cdot \varepsilon \cdot
|y_{t}^{'}(t_0,x,\varepsilon)|+
$$
\begin{equation}
\label{kop_9}
 + K \cdot \max\limits_{(t,x \in G)} |K_{1}(t,0,x) \cdot
\Delta_{0}(x) - \Delta(t,x)| + K \cdot
\left(1+|y_{t}^{'}(t_0,x,\varepsilon)| \right) \cdot
e^{-\frac{\gamma}{\varepsilon}(t-t_0)},
\end{equation}
where $K$ - some constant independent on $t$ and $\varepsilon$,
$G_1=\{(t,x): 0 <t_0 \leq t \leq 1, \lambda \leq x \leq \lambda +
1\}$.
\end{thm}
{\it Proof}. Indeed, in (\ref{kop_1}) assign
$y(t,x,\varepsilon)=y_0(t,x)+u(t,x,\varepsilon)$, and taking into
consideration (\ref{kop_7}), (\ref{kop_8}), we obtain for
$u(t,x,\varepsilon)$ following problem
\begin{equation}
\label{kop_10} L_{\varepsilon}u=f(t,x,\varepsilon)+
\int_{t_0}^{t}\left(K_{1}(t,s,x)H[u(s,x,\varepsilon)]+K_0(t,s,x)u(s,x,\varepsilon)\right)ds,
\end{equation}
\begin{equation}
\label{kop_11} u(t_0,x,\varepsilon)=0, \qquad
u_{t}^{'}(t_0,x,\varepsilon)=y_{t}^{'}(t_0,x,\varepsilon),
\end{equation}
where function $f(t,x,\varepsilon)$ has a representation
$$
f(t,x,\varepsilon)=K_1(t,0,\varphi) \cdot \Delta_0(\psi)+
\int_{t_0}^{t}\left(K_{0}(t,s,\varphi)-\frac{\partial
K_1(t,s,\varphi)}{\partial s}\right)u(s,\varphi,\varepsilon)ds -
\Delta(t,x) - \varepsilon H^{2}[y_0(t,x)],
$$
and estimation
\begin{equation}
\label{kop_12} |f(t,x,\varepsilon)| \leq K \cdot \varepsilon \cdot
|\ln \varepsilon|+ \max\limits_{(t,x \in G)} |K_{1}(t,0,x) \cdot
\Delta_{0}(x) - \Delta(t,x)|, \quad (t,x) \in G,
\end{equation}
where $\varphi=\varphi(t,\psi)$ - a solution of characteristic
equation (\ref{kop_3}), $\psi=\psi(t,x)=x_0$ - first integral of
equation (\ref{kop_1}) \cite{Tazhim}. Applying to the problem
(\ref{kop_10}), (\ref{kop_11}) theorem 1 from \cite{Tazhim}, and
taking into consideration (\ref{kop_12}), we obtain (9). Theorem
is proved.

\section{Initial jumps of solutions and integral term}
The aim of this section is to define the conditions in the
presence of which the solution of perturbed Cauchy problem goes to
the solution of unperturbed problem.

Taking into consideration (\ref{kop_9}), assign
\begin{equation}
\label{kop_13} \Delta(t,x)=\Delta_0(x) \cdot K_1(t,0,x).
\end{equation}
Then from theorem1, we obtain, that
$$
\lim\limits_{x \to 0} y(t,x,\varepsilon) = y_0(t,x), \quad
\lim\limits_{x \to 0} y_{t}^{'}(t,x,\varepsilon) = \frac{\partial
y_0(t,x)}{\partial t}, \quad (t,x) \in G_1.
$$
Further, for define $\Delta_0(x)$, integrate equation
(\ref{kop_1}) along characteristic $x=\varphi(t,\psi)$ on $t$ from
$0$ to $t_0$. Then we have
$$
\varepsilon \cdot H[y_0(t_0,\varphi,\varepsilon)] - \pi_{1}(\psi)
+
A(t_0,\varphi)y_0(t_0,\varphi,\varepsilon)-A(0,\varphi)\pi_{0}(\psi)-
\int_{0}^{t_0}\left(A_{t}^{'}(t,\varphi)-B(t,\varphi)\right)\times
$$
\begin{equation}
\label{kop_14} \times y(t,\varphi,\varepsilon)dt =
\int_{0}^{t_0}\left( F(t,\varphi)+\int_{0}^{t}\left(
K_1(t,s,\varphi)H[y(s,\varphi,\varepsilon)]+K_0(t,s,\varphi)y(s,\varphi,\varepsilon)
\right)ds\right)dt
\end{equation}
From (\ref{kop_14}), passage to the limit when $\varepsilon \to
0$, and taking into consideration (\ref{kop_7}), ( \ref{kop_9}),
(\ref{kop_13}), obtain
\begin{equation}
\label{kop_15} \Delta_0(x)=\frac{\pi_1(\psi)}{A(0,\psi)}, \quad
\Delta(t,x)=\frac{\pi_1(\psi)}{A(0,\psi)} \cdot K_1(t,0,\varphi).
\end{equation}
Thus, in $G_1$, if equalities (\ref{kop_15}) is satisfied, then
difference between solution $y(t,x,\varepsilon)$ of problem
(\ref{kop_1}), (\ref{kop_2}) and solution $y_0(t,x)$ of problem
(\ref{kop_6}), (\ref{kop_7}) will be enough small with
$\varepsilon$.

\bibliographystyle{plain}
\bibliography{nonexistence}
\end{document}